\documentclass[reqno]{amsart}
\usepackage{amsmath, amsthm, amssymb, amscd}

\usepackage[a4paper,width=150mm,top=25mm,bottom=25mm,centering]{geometry}

\usepackage{enumitem}
\usepackage{graphicx} % Required for inserting images
\usepackage[dvipsnames]{xcolor}
\usepackage{comment}
\usepackage{cases}
\usepackage{array}
\usepackage{multirow}
\usepackage{microtype}
\usepackage{longtable}

\usepackage{tikz}
\usepackage{pgf}
\usepackage{tikzscale}
\usepackage{tikz-cd}
\usetikzlibrary{decorations.pathreplacing,decorations.markings}
\usetikzlibrary{arrows, knots, external, positioning}

\makeatletter
\tikzset{column sep/.code=\def\pgfmatrixcolumnsep{\pgf@matrix@xscale*(#1)},
  row sep/.code   =\def\pgfmatrixrowsep{\pgf@matrix@yscale*(#1)},
  matrix xscale/.code=%
    \pgfmathsetmacro\pgf@matrix@xscale{\pgf@matrix@xscale*(#1)},
  matrix yscale/.code=%
    \pgfmathsetmacro\pgf@matrix@yscale{\pgf@matrix@yscale*(#1)},
  matrix scale/.style={/tikz/matrix xscale={#1},/tikz/matrix yscale={#1}}}
\def\pgf@matrix@xscale{1}
\def\pgf@matrix@yscale{1}
\makeatother

\usepackage{url}
\usepackage{caption} % to give a padding between a table and its caption

\newtheorem{theorem}{Theorem}[section]
\newtheorem{conjecture}[theorem]{Conjecture}
\newtheorem{lemma}[theorem]{Lemma}
\newtheorem{proposition}[theorem]{Proposition}
\newtheorem{corollary}[theorem]{Corollary}
\newtheorem{question}[theorem]{Question}

\theoremstyle{definition}

\newtheorem{remark}[theorem]{Remark}

\numberwithin{equation}{section}

\usepackage[bookmarksnumbered=true, colorlinks=true]{hyperref}
\AtBeginDocument{
    \hypersetup{
        citecolor=BlueGreen,
        linkcolor=BlueGreen,
        urlcolor=BlueGreen
    }
}

\begin{document}

\title[Topology of symplectic 4-orbifolds with vanishing canonical class]{Topology of symplectic 4-orbifolds with vanishing canonical class}

\author{Woohyeok Jo}
\address{HUN-REN R\'enyi Institute of Mathematics, H-1053 Budapest, Re\'altanoda utca 13--15, Hungary}
\email{jo.woohyeok@renyi.hu}

\author{Jongil Park}
\address{Department of Mathematical Sciences and Research Institute of Mathematics, Seoul National University, Seoul 08826, Republic of Korea}
\email{jipark@snu.ac.kr}

\thanks{}
% 57K43  	Symplectic structures in 4 dimensions
% 57R18     Topology and geometry of orbifolds
% 57R20     Characteristic classes and numbers in differential topology
% 14J28     K3 surfaces and Enriques surfaces

\date\today

\begin{abstract}
We prove that every closed, effective symplectic $4$-orbifold with isolated singularities has nonnegative orbifold Euler characteristic $e_{\mathrm{orb}}$ whenever its canonical class vanishes in rational cohomology. We also establish uniform bounds on the global index, the orders of the isotropy groups, the number of singular points, and the second Betti number. Consequently, only finitely many singularity types occur in this class of orbifolds. In the equality case $e_{\mathrm{orb}}=0$, we further prove that the orbifold admits a connected finite orbifold covering by a closed smooth symplectic $4$-manifold with torsion canonical class and zero Euler characteristic.
\end{abstract}

\maketitle

\section{Introduction}

For a compact K\"ahler surface $S$ of nonnegative Kodaira dimension, the Bogomolov--Miyaoka--Yau \textup{(BMY)} inequality and the classification of complex surfaces give \begin{equation}\label{eq:complex_BMY}
c_1^2(S)\leq3c_2(S), \quad c_2(S)\geq0;
\end{equation}
see \cite{Miyaoka-1977,Yau-1977,Barth-Hulek-Peters-VanDeVen-2004}. Note that $c_1^2(S)$ equals $K_S^2$, where $K_S$ is the canonical class of $S$, and that $c_2(S)$ equals the topological Euler characteristic $e(S)$.

For surfaces with quotient singularities, the corresponding inequalities involve the orbifold second Chern number $c_2^{\mathrm{orb}}(S)$, or equivalently the orbifold Euler characteristic $e_{\mathrm{orb}}(S)$. If $S$ is a normal projective complex surface with quotient singularities and $\Gamma_p$ is the isotropy group at $p\in\operatorname{Sing}(S)$, then \begin{equation}\label{eq:orbifold_euler}
e_{\mathrm{orb}}(S)=e(S)-\sum_{p\in\operatorname{Sing}(S)} \left(1-\frac1{|\Gamma_p|}\right),
\end{equation}
where $e(S)$ is the ordinary Euler characteristic of $S$ \cite{Satake-1957, Kollar-2008}. We recall the following two forms of the orbifold BMY inequality.

\begin{theorem}\label{thm:orbifold_BMY_inequalities}
Let $S$ be a normal projective complex surface with quotient singularities. 
\begin{enumerate}[label=\textup{(\alph*)}]
\item \textup{(Orbifold BMY inequality \cite{Sakai-1980,Miyaoka-1984,KoNS-1989,Megyesi-1999})} If $K_S$ is nef, then \[
K_S^2\leq3e_{\mathrm{orb}}(S).
\]
In particular, $e_{\mathrm{orb}}(S)\geq0$.
\item \textup{(Weak orbifold BMY inequality \cite{KeM-1999})} If $-K_S$ is nef, then $e_{\mathrm{orb}}(S)\geq0$.
\end{enumerate}
\end{theorem}

These inequalities are particularly useful for algebraic $\mathbb{Q}$-homology projective planes, namely, normal projective complex surfaces with the same rational homology as $\mathbb{CP}^2$. For such a surface $S$ with quotient singularities, either $K_S$ or $-K_S$ is nef, and Theorem \ref{thm:orbifold_BMY_inequalities} therefore gives $e_{\mathrm{orb}}(S)\geq0$. The role of these inequalities in the study of such surfaces and related problems is discussed in Koll\'ar's survey \cite{Kollar-2008}. The orbifold BMY inequality is also an important ingredient in the resolution of the algebraic Montgomery--Yang problem \cite{Jo-Park-Park-2026}.

For symplectic $4$-manifolds of nonnegative symplectic Kodaira dimension, the inequalities \eqref{eq:complex_BMY} are conjectured to hold:
\begin{conjecture}[{\cite[Problem 4.90]{K3-2026}}]\label{conj:smooth_chern_inequalities}
Let $M$ be a closed symplectic $4$-manifold of nonnegative symplectic Kodaira dimension. Then \[
c_1^2(M)\leq3c_2(M),\quad c_2(M)\geq0.
\]
\end{conjecture}

We consider analogous questions for symplectic $4$-orbifolds with isolated singularities. Note that every normal projective complex surface with quotient singularities is an example after choosing an orbifold K\"ahler form. For such a symplectic $4$-orbifold $(X,\omega)$, 
we write \[
K_X=-c_1^{\mathrm{orb}}(TX,J)\in H^2(X;\mathbb{Q})
\]
for the rational canonical class, where $J$ is any $\omega$-compatible orbifold almost complex structure. Formula \eqref{eq:orbifold_euler} also defines $e_{\mathrm{orb}}(X)$, and orbifold integration gives $c_1^{\mathrm{orb}}(X)^2=K_X^2$ and $c_2^{\mathrm{orb}}(X)=e_{\mathrm{orb}}(X)$. For normal projective surfaces, nefness of $K_X$ implies $K_X^2\geq0$ and $K_X\cdot[\omega]\geq0$ for every orbifold K\"ahler form $\omega$, while nefness of $-K_X$ implies $K_X^2\geq0$ and $K_X\cdot[\omega]\leq0$. These conditions motivate the following analogues of the two parts of Theorem \ref{thm:orbifold_BMY_inequalities}.

\begin{question}\label{q:numerical_orbifold_bmy}
Let $(X,\omega)$ be a closed, connected symplectic $4$-orbifold with isolated singularities. If $K_X^2\geq0$ and $K_X\cdot[\omega]\geq0$, then does $K_X^2\leq 3e_{\mathrm{orb}}(X)$ hold?
\end{question}

\begin{question}\label{q:weak_orbifold_bmy_negative}
Let $(X,\omega)$ be a closed, connected symplectic $4$-orbifold with isolated singularities. If $K_X^2\geq0$ and $K_X\cdot[\omega]\leq0$, then does $e_{\mathrm{orb}}(X)\geq0$ hold?
\end{question}

In this paper, we study the case $K_X=0\in H^2(X;\mathbb{Q})$. For a symplectic $4$-manifold, the vanishing of the rational canonical class is equivalent to saying that the integral canonical class is torsion. Such a manifold is symplectically minimal of Kodaira dimension zero, and we have the following result of Li.

\begin{theorem}[{\cite[Corollary 1.2]{Li-Quaternionic-2006}}]\label{thm:li_smooth_nonnegativity}
Let $M$ be a closed symplectic $4$-manifold. If its canonical class $K_M\in H^2(M;\mathbb{Z})$ is torsion, then $e(M)\geq0$.
\end{theorem}

Our first main result extends this nonnegativity to symplectic $4$-orbifolds.

\begin{theorem}\label{thm:main}
Let $X$ be a closed symplectic $4$-orbifold with isolated singularities. If $K_X=0\in H^2(X;\mathbb{Q})$, then \[
e_{\mathrm{orb}}(X)\geq0.
\]
Equivalently, $c_2^{\mathrm{orb}}(X)\geq0$.
\end{theorem}

Only the rational canonical class is required to vanish; the canonical orbifold line bundle need not be trivial. Since $K_X^2=0$, Theorem \ref{thm:main} gives affirmative answers to Questions \ref{q:numerical_orbifold_bmy} and \ref{q:weak_orbifold_bmy_negative} in this case.

Our next result concerns the possible singularities of $X$. Each isolated singularity $p\in\operatorname{Sing}(X)$ has a local model $\mathbb{C}^2/\Gamma_p$, where $\Gamma_p\subset\mathrm{U}(2)$ is a finite group acting freely away from the origin; see Section \ref{sec:preliminaries}. The \emph{local index} at $p$ is defined by \[
r_p:=[\Gamma_p:\Gamma_p\cap\mathrm{SU}(2)].
\]
The \emph{global index} $r(X)$ is defined to be the least common multiple of the local indices, with $r(X)=1$ when $X$ is smooth. Our second main result is the following finiteness theorem.

\begin{theorem}\label{thm:finiteness}
There exists a universal positive integer $B$ such that, for every closed, connected symplectic $4$-orbifold $X$ with isolated singularities and $K_X=0\in H^2(X;\mathbb{Q})$,
\[
r(X)\mid B,\quad b_2(X)\leq22,\quad |\operatorname{Sing}(X)|\leq48,
\]
and $|\Gamma_p|\leq120B$ for every $p\in\operatorname{Sing}(X)$. Consequently, only finitely many singularity types occur in this class.
\end{theorem}

By the \emph{singularity type} of $X$, we mean the unordered collection, counted with multiplicity, of the types of its local quotient singularities. An explicit choice of $B$ is given in Remark \ref{rmk:choice_of_B}. The bound for $b_2(X)$ is sharp, since $b_2=22$ for smooth K3 surfaces. We do not claim optimality for the bounds $r(X)\leq B$ and $|\operatorname{Sing}(X)|\leq 48$. 

\vspace{2mm}
A \emph{symplectic $\mathbb{Q}$-homology projective plane} is a closed symplectic $4$-orbifold $X$ with isolated singularities and $H_*(X;\mathbb{Q})\cong H_*(\mathbb{CP}^2;\mathbb{Q})$. For this class, our results have the following consequence.

\begin{corollary}\label{cor:canonical_zero_qhpp}
Let $X$ be a symplectic $\mathbb{Q}$-homology projective plane with $K_X=0\in H^2(X;\mathbb Q)$. Then
\[
e_{\mathrm{orb}}(X)\geq0,\quad |\operatorname{Sing}(X)|\leq5.
\]
Both bounds are sharp. Moreover, only finitely many singularity types occur among such orbifolds.
\end{corollary}

In \cite{Jo-Park-2026}, these inequalities are proved for symplectic $\mathbb{Q}$-homology projective planes with only rational double points, without any hypothesis on the canonical class. Corollary \ref{cor:canonical_zero_qhpp} allows arbitrary quotient singularities under the assumption that the rational canonical class vanishes.

\vspace{2mm}
Our third main result concerns the equality case $e_{\mathrm{orb}}=0$. In the compact K\"ahler case, a normal surface $S$ with quotient singularities and $K_S=0\in H^2(S;\mathbb{Q})$ satisfies $e_{\mathrm{orb}}(S)=0$ if and only if it is a finite quotient of a complex $2$-torus by a holomorphic action free outside finitely many points \cite{Kobayashi-Todorov-1987}. In particular, it admits a finite holomorphic orbifold covering by a smooth complex surface. We establish a finite smooth covering statement in the symplectic setting.

\begin{theorem}\label{thm:zero_euler_smooth_cover}
Let $X$ be a closed, connected symplectic $4$-orbifold with isolated singularities. If $K_X=0\in H^2(X;\mathbb{Q})$ and $e_{\mathrm{orb}}(X)=0$, then there exists a finite-degree orbifold covering $M \to X$, where $M$ is a closed, connected symplectic $4$-manifold such that $e(M)=0$ and $K_M \in H^2(M;\mathbb{Z})$ is torsion.
\end{theorem}

The conclusion need not hold when $e_{\mathrm{orb}}(X)>0$ even for projective surfaces; see Remark \ref{rmk:positive_euler_no_smooth_cover}.

\subsection*{Outline of the proofs}
For Theorem \ref{thm:main}, Chen's Calabi--Yau cover reduces the problem to the case where $X$ has only rational double points. If the orbifold Euler characteristic were negative, Bauer's bound would imply a uniform bound on the degrees of finite orbifold covers. The smooth locus of a cover of maximal degree has vanishing first integral homology. The minimal resolution then has the intersection lattice of a K3 surface, and its exceptional lattice is primitive. Shimada's criterion then realizes the same singularity type on a normal projective K3 surface, and the algebraic orbifold BMY inequality gives a contradiction.

For Theorem \ref{thm:finiteness}, we treat the case $r(X)=1$ separately. When $r(X)>1$, we apply Chen--Kwasik's rigidity theorem to the minimal resolution of the Calabi--Yau cover. We then use the fact that there is a uniform bound for the orders of finite-order matrices in $\operatorname{GL}(22,\mathbb{Z})$ to obtain a uniform bound for $r(X)$. Combining this with bounds on the exceptional rank gives a uniform bound on the orders of the isotropy groups. Finally, Theorem \ref{thm:main} gives a uniform bound on the number of singular points. Together, these bounds imply that only finitely many singularity types occur.

For Theorem \ref{thm:zero_euler_smooth_cover}, we use Kobayashi--Todorov's characterization to show that the smooth locus of a normal projective K3 surface with zero orbifold Euler characteristic has nonvanishing first integral homology. From this, we deduce that the degrees of finite covers are unbounded in the equality case $e_{\mathrm{orb}}=0$. On the other hand, Bauer's bound limits the total exceptional rank in the resolution of every finite cover. A regular cover of a sufficiently large degree then provides the required smooth cover.

\vspace{2mm}
We conclude the introduction by posing further questions.
First, Theorem \ref{thm:finiteness} naturally leads to the following classification problem.

\begin{question}\label{q:singularity_classification}
Which singularity types occur on closed, connected symplectic $4$-orbifolds with isolated singularities and vanishing rational canonical class?
\end{question}

Related questions concerning singular sets and their realizations have also been posed in \cite[Problem 1.3]{Chen-2020}. In addition, any such orbifold $X$ satisfies homological restrictions: $1\leq b_2^+(X)\leq3$, with $b_2^+(X)=1$ whenever $X$ has a singularity that is not a rational double point; see Proposition \ref{prop:orbifold_bplus}.

On the other hand, complex surfaces with numerically trivial canonical divisor have been extensively studied \cite{Zhang-1991,Zhang-1993,Blache-1995,Oguiso-Zhang-1996,Oguiso-Zhang-1999,Kudryavtsev-2002,Wang-2011,Esser-2022,Schutt-2023}. Notably, Blache \cite[Theorem C]{Blache-1995} proved that a minimal normal compact complex Moishezon surface $S$ with only log canonical singularities and Kodaira dimension zero has Cartier index \[
I(S):=\min\{m>0:mK_S\text{ is Cartier}\}\leq21.
\]
For surfaces with quotient singularities, $I(S)$ is the least common multiple of the local indices and hence agrees with the global index $r(S)$ as defined above. Later, Esser--Totaro--Wang \cite[Proposition 6.1]{Esser-2022} proved that every normal projective surface $S$ with only quotient singularities such that $K_S$ is $\mathbb{Q}$-linearly equivalent to zero satisfies \[
m(S):=\min \{m>0:mK_S\text{ is linearly equivalent to zero}\} \leq 19.
\]
Since $r(S)=I(S)$ divides $m(S)$, this gives $r(S)\leq 19$. Moreover, this bound is sharp: examples with $r(S)=19$ are given in \cite{Zhang-1991, Blache-1995}. This motivates the following question. 

\begin{question}
Does every closed, connected symplectic $4$-orbifold $X$ with isolated singularities and $K_X=0\in H^2(X;\mathbb{Q})$ satisfy $r(X)\leq 19$?
\end{question}

\vspace{2mm}
The organization of this paper is as follows. 
In Section \ref{sec:preliminaries}, we fix conventions and collect the results used in the proofs. In Sections \ref{sec:vanishing_canonical} and \ref{sec:finiteness}, we prove Theorems \ref{thm:main} and \ref{thm:finiteness}, respectively. We also prove Corollary \ref{cor:canonical_zero_qhpp} in the latter section. Finally, we prove Theorem \ref{thm:zero_euler_smooth_cover} in Section \ref{sec:equality_case}.

\subsection*{Acknowledgments}
The authors thank all members of the 4-manifold topology group at Seoul National University (SNU) for their comments throughout this work. Woohyeok Jo was supported by ERC Advanced Grant KnotSurf4d. Jongil Park was supported by the National Research Foundation of Korea (NRF) grant funded by the Korean government (RS-2024-00392067). He is also affiliated with the Research Institute of Mathematics at SNU.

\section{Preliminaries}\label{sec:preliminaries}

In this section, we review the background materials and collect the results used in the proofs.

\subsection{Lattices}\label{subsec:lattices}
We first fix our conventions concerning lattices and record a criterion for primitivity.

A \emph{lattice} is a pair $(L,Q)$ consisting of a free abelian group $L$ of finite rank and a nondegenerate symmetric bilinear form $Q:L\times L\to\mathbb Q$. For $v,w\in L$, we write $v\cdot w:=Q(v,w)$ and $v^2:=Q(v,v)$. An \emph{embedding} (respectively, an \emph{isomorphism}) of lattices is an injective group homomorphism (respectively, a group isomorphism) that preserves the bilinear forms. The \emph{signature} of a lattice is the pair of positive and negative indices of the induced form on $L\otimes\mathbb{R}$.

A lattice $(L,Q)$ is \emph{integral} if $Q$ takes values in $\mathbb{Z}$, and an integral lattice $(L,Q)$ is \emph{even} if $v^2\in2\mathbb{Z}$ for every $v\in L$. For an integral lattice $(L,Q)$, put $L_{\mathbb{Q}}:=L\otimes\mathbb{Q}$ and extend $Q$ to a $\mathbb{Q}$-bilinear form $L_{\mathbb{Q}}\times L_{\mathbb{Q}}\to \mathbb{Q}$. The \emph{dual lattice} of $(L,Q)$ is
\[
L^*:=\{x\in L_{\mathbb{Q}}:x\cdot v\in\mathbb{Z}\text{ for every }v\in L\},
\]
with the restricted bilinear form. The map $x\mapsto Q(x,-)$ canonically identifies the underlying abelian group of $L^*$ with $\operatorname{Hom}(L,\mathbb{Z})$. We have $L\subset L^*$, and $L$ is \emph{unimodular} if $L=L^*$, equivalently if $\det Q=\pm1$, where $\det Q$ is the determinant of the Gram matrix with respect to any basis for $L$. 

A sublattice $L\subset\Lambda$ is \emph{primitive} if $\Lambda/L$ is torsion-free, or equivalently if $(L\otimes\mathbb{Q})\cap\Lambda=L$, where the intersection is taken in $\Lambda\otimes\mathbb{Q}$. The lattice $(L\otimes\mathbb{Q})\cap\Lambda$ is called the \emph{primitive closure} of $L$ in $\Lambda$.

\begin{lemma}\label{lem:primitive_restriction}
Let $\Lambda$ be a unimodular integral lattice and let $L\subset\Lambda$ be a (nondegenerate) sublattice. Then $L$ is primitive in $\Lambda$ if and only if the restriction map \[
r_L:\Lambda\longrightarrow L^*,\quad r_L(x)(v)=x\cdot v\quad(v\in L),
\]
is surjective, where $L^*$ is identified with $\operatorname{Hom}(L,\mathbb{Z})$.
\end{lemma}
\begin{proof}
Unimodularity identifies $\Lambda$ with $\operatorname{Hom}(\Lambda,\mathbb{Z})$, under which $r_L$ is the restriction map induced by $L\hookrightarrow \Lambda$. Applying $\operatorname{Hom}(-,\mathbb{Z})$ to $0\to L\to\Lambda\to\Lambda/L\to0$ gives \[
\operatorname{coker}(r_L)\cong\operatorname{Ext}^1_{\mathbb{Z}}(\Lambda/L,\mathbb{Z}).
\]
Since $\Lambda/L$ is finitely generated, this group vanishes if and only if $\Lambda/L$ is torsion-free.
\end{proof}

We will use the following lattices and conventions.
\begin{itemize}[leftmargin=2em]
    \item For a compact, oriented $4$-manifold $M$ with $H^1(\partial M;\mathbb{Z})=0$, we denote its intersection lattice on $H_2(M;\mathbb{Z})/\text{tor}$ by $Q_M$. The universal coefficient theorem identifies $H^2(M;\mathbb{Z})/\text{tor}$ with $\operatorname{Hom}(H_2(M;\mathbb{Z})/\text{tor},\mathbb{Z})$ and hence with $Q_M^*$. If $M$ is closed, then $Q_M$ is unimodular. We write $b_2^+(M)$ and $b_2^-(M)$ for its positive and negative indices, respectively, and we write $\sigma(M)=b_2^+(M)-b_2^-(M)$ for the signature of $M$.
    \item We use $A_n$ $(n\geq1)$, $D_n$ $(n\geq4)$, and $E_6,E_7,E_8$ for the negative definite root lattices of the corresponding Dynkin types. An \emph{\textup{ADE} root lattice} is an orthogonal direct sum of these lattices. A vector of square $-2$ is called a \emph{root}.
    \item We denote by $U$ the rank-two hyperbolic lattice with Gram matrix \[
    \begin{pmatrix}
        0 & 1 \\ 1 & 0
    \end{pmatrix}.
    \]
\item The \emph{\textup{K3} lattice} is \[
\Lambda_{\mathrm{K3}}:=3U\oplus 2E_8=U^{\oplus3}\oplus E_8^{\oplus2}.
\]
With our negative definite convention for $E_8$, this is the even unimodular lattice of signature $(3,19)$, unique up to isomorphism. It is isomorphic to the intersection lattice of every complex $\textup{K3}$ surface. 
\end{itemize}

\subsection{Symplectic 4-orbifolds}
\label{subsec:symplectic_orbifolds}

Throughout this paper, we assume that all orbifolds are effective, and homology and cohomology groups are understood to be those of their underlying topological spaces unless otherwise specified. The isotropy group at a point $p$ is denoted by $\Gamma_p$. The points in an orbifold $X$ with trivial isotropy group form an open dense subset, called the \emph{smooth locus} and denoted by $X^\circ$. 

Let $(X,\omega)$ be a symplectic $4$-orbifold with only isolated singularities. For $p\in \operatorname{Sing}(X)$, the equivariant Darboux theorem and a $\Gamma_p$-invariant compatible linear complex structure give a local model: \[
(B_p/\Gamma_p,\omega_{\mathrm{std}}), \quad \Gamma_p\subset\mathrm{U}(2),
\]
where $B_p\subset\mathbb{C}^2$ is a ball centered at the origin and $\Gamma_p$ is finite and acts freely on $\mathbb{C}^2-\{0\}$. Thus, with the induced local complex structure, the germ of $X$ at $p$ is a surface quotient singularity $(\mathbb{C}^2/\Gamma_p,0)$ with link $S^3/\Gamma_p$. In particular, the isotropy group $\Gamma_p$ is isomorphic to the \emph{local fundamental group}, i.e., the fundamental group of the link. 

The point $p\in \operatorname{Sing}(X)$ is called a \emph{rational double point} if $\Gamma_p\subset\mathrm{SU}(2)$. The exceptional set of its minimal algebraic resolution consists of smooth rational $(-2)$-curves with dual graph of type $A_n$, $D_n$, $E_6$, $E_7$, or $E_8$. The associated intersection lattice is the negative definite root lattice of the same type. We recall the orders of the corresponding isotropy groups:
\begin{equation}\label{eq:ade_group_orders}
\begin{array}{c|ccccc}
\text{Type}&A_n&D_n&E_6&E_7&E_8\\ \hline
|\Gamma_p|&n+1&4(n-2)&24&48&120
\end{array}
\end{equation}

For a closed symplectic $4$-orbifold $(X,\omega)$ with isolated singularities, a \emph{symplectic resolution} consists of a closed symplectic $4$-manifold $(M,\omega_M)$ and a continuous map $\rho:M\to X$ restricting to a diffeomorphism $M-E\to X^\circ$, where $E:=\rho^{-1}(\operatorname{Sing}(X))$ is a union of symplectic surfaces, such that $\omega_M=\rho^*\omega$ outside a small neighborhood of $E$. Such resolutions exist \cite[Theorem 3.3]{Cavalcanti-Fernandez-Munoz}. Near each singular point, we can choose the resolution to be the minimal algebraic resolution of the corresponding quotient singularity. We call the resulting resolution the \emph{minimal symplectic resolution}; see also \cite[Section 2.2]{Jo-Park-2026}. The exceptional set over each singular point thus has the weighted dual graph of its minimal algebraic resolution.

Since each link is a rational homology $3$-sphere, the orbifold $X$ is a rational homology $4$-manifold. In particular, its ordinary cup-product pairing on $H^2(X;\mathbb{Q})$ is nondegenerate. If $\rho:M\to X$ is the minimal symplectic resolution, the pullback embeds $H^2(X;\mathbb{Q})$ into $H^2(M;\mathbb{Q})$, and its orthogonal complement is the negative definite subspace spanned by the Poincar\'e duals of the exceptional curves. Thus $b_2(M)=b_2(X)+L$ and $b_2^+(M)=b_2^+(X)$, where $L$ is the total number of exceptional curves and $b_2^+(X)$ is the positive index of the real cup-product pairing.

Choose an $\omega$-compatible orbifold almost complex structure $J$. We denote the orbifold first Chern class of its canonical orbifold line bundle $\Lambda_{\mathbb{C}}^{2}(T^{*}X,J)$ by $K_X\in H^2(X;\mathbb{Q})$ and call it the \emph{canonical class} of $(X,\omega)$. This class is independent of the choice of $J$. For a symplectic $4$-manifold $M$, we retain the usual integral canonical class $K_M=-c_1(TM,J)\in H^2(M;\mathbb{Z})$ and use the same notation for its image in $H^2(M;\mathbb{Q})$.

\begin{lemma}\label{lem:rdp_canonical_pullback}
Let $\rho:M\to X$ be the minimal symplectic resolution of a closed symplectic $4$-orbifold $X$ with only rational double points. Then \[
K_M=\rho^*K_X \in H^2(M;\mathbb{Q}).
\]
In particular, if $K_X=0\in H^2(X;\mathbb{Q})$, then $K_M\in H^2(M;\mathbb{Z})$ is torsion.
\end{lemma}
\begin{proof}
This follows from the fact that $K_M$ pairs trivially with every exceptional sphere by the adjunction formula, since the difference $K_M-\rho^*K_X$ is a $\mathbb{Q}$-linear combination of the Poincar\'e duals of the exceptional spheres. 
\end{proof}

For a finite symplectic orbifold covering $q:Y\to X$ of degree $d$, we have $K_Y=q^*K_X$. When $X$ and $Y$ are closed, the underlying map has degree $d$, and hence \begin{equation}\label{eq:cover_pairing}
(q^*u)\cdot(q^*v)=d(u\cdot v) \quad (u,v\in H^2(X;\mathbb{Q})).
\end{equation}
Since the pairings are nondegenerate, it follows that $q^*:H^2(X;\mathbb{Q})\to H^2(Y;\mathbb{Q})$ is injective. We also have the multiplicativity $e_{\mathrm{orb}}(Y)=d\cdot e_{\mathrm{orb}}(X)$.

\vspace{2mm}
Let $(X,\omega)$ be a closed symplectic $4$-orbifold with isolated singularities. For $p\in\operatorname{Sing}(X)$, set \[
H_p:=\Gamma_p\cap\mathrm{SU}(2),\quad r_p:=[\Gamma_p:H_p]=|\det(\Gamma_p)|.
\]
The integer $r_p$ is called the \emph{local index} at $p$. It equals one precisely at rational double points. The \emph{global index} of $X$ is defined by
\begin{equation}\label{eq:local_index_lcm}
r(X):=\underset{p\in \operatorname{Sing}(X)}{\operatorname{lcm}}r_p,
\end{equation}
with $r(X)=1$ when $X$ is smooth. Thus $r(X)=1$ precisely when every singularity is a rational double point, including the smooth case. We now recall Chen's Calabi--Yau covering theorem together with its properties.

\begin{theorem}[{\cite{Chen-2020}}]\label{thm:chen_canonical_cover}
Let $(X,\omega)$ be a closed, connected symplectic $4$-orbifold with isolated singularities and $K_X=0\in H^2(X;\mathbb{Q})$. Suppose that $r:=r(X)>1$. Then $b_1(X)=0$, and there exists a connected cyclic symplectic orbifold covering $q:Y\to X$ of degree $r$ such that $Y$ has only rational double points and a trivial canonical orbifold line bundle. Moreover:
\begin{enumerate}[label=\textup{(\alph*)}, leftmargin=2em]
\item For each $p\in\operatorname{Sing}(X)$, there are $r/r_p$ points over $p$. Near each such point, the covering is modeled on \[
B_p/H_p\longrightarrow B_p/\Gamma_p.
\]
In particular, each point of $q^{-1}(p)$ has isotropy group $H_p$, and is smooth if $H_p=\{1\}$.
\item The minimal symplectic resolution $M\to Y$ can be chosen equivariantly, with $K_M=0\in H^2(M;\mathbb{Z})$, so that the deck transformations lift to an effective symplectic $\mathbb{Z}_r$-action on $M$. The resolution $M$ has the integral homology of a $\operatorname{K3}$ surface, except possibly in the following two cases: \begin{enumerate}[label=\textup{(\roman*)}]
\item $X$ has nine singularities of type $\frac{1}{3}(1,1)$ and $r=3$;
\item $X$ has five singularities of type $\frac{1}{5}(1,2)$ and $r=5$.
\end{enumerate}
In either exceptional case, $M=Y$ is smooth and has the rational homology of $T^4$.
\end{enumerate}
\end{theorem}

The orbifold $Y$ is called the \emph{Calabi--Yau cover} of $X$. 

\vspace{2mm}
We will also use the following result of Bauer in the proofs below. Its proof uses the $\operatorname{Pin}(2)$-equivariant stable cohomotopy refinement of the Seiberg--Witten monopole map and equivariant obstruction theory, together with Taubes's theorem. See also \cite[Corollary 1.2]{Li-Quaternionic-2006}.

\begin{theorem}[{\cite[Corollary 1.2]{Bauer-2008}}]\label{thm:Bauer}
Let $M$ be a closed, connected, symplectic $4$-manifold. If its canonical class $K_M\in H^2(M;\mathbb{Z})$ is torsion, then $b_2^+(M)\leq3$.
\end{theorem}

We next recall a homological rigidity theorem for finite symplectic group actions. An action is \emph{homologically trivial over $\mathbb Q$} if every group element induces the identity on $H_*(-;\mathbb{Q})$.

\begin{theorem}[{\cite[Theorem~A]{Chen-Kwasik-2007}}]\label{thm:chen_kwasik_rigidity}
Let $(M,\omega)$ be a closed, connected symplectic $4$-manifold with $K_M=0\in H^2(M;\mathbb Q)$, $\sigma(M)\neq0$, and $b_2^+(M)\geq2$. Then every finite symplectic group action on $M$ that is homologically trivial over $\mathbb{Q}$ must be trivial.
\end{theorem}

We end this subsection by noting the following consequence of Bauer's and Chen's results.

\begin{proposition}\label{prop:orbifold_bplus}
Let $(X,\omega)$ be a closed, connected symplectic $4$-orbifold with isolated singularities and $K_X=0\in H^2(X;\mathbb{Q})$. Then $1\leq b_2^+(X)\leq3$. Furthermore, if there is a singularity that is not a rational double point, then $b_2^+(X)=1$.
\end{proposition}
\begin{proof}
Since $[\omega]^2>0$, we have $b_2^+(X)\geq1$. Let $\rho:M\to X$ be the minimal symplectic resolution. Its exceptional lattice is negative definite, so $b_2^+(M)=b_2^+(X)$. If every singularity is a rational double point, Lemma \ref{lem:rdp_canonical_pullback} shows that $K_M$ is torsion, and Theorem \ref{thm:Bauer} gives $b_2^+(X)\leq3$. If there is a singularity that is not a rational double point, then $M$ is rational or ruled by \cite[Proposition 1.0]{Chen-2020}, so $b_2^+(X)=b_2^+(M)=1$.
\end{proof}

\subsection{K3 surfaces}
\label{subsec:k3_configurations}

In this subsection we recall Shimada's criterion \cite{Shimada-2007} using notation adapted to the present paper. A \emph{normal projective $\textup{K3}$ surface} means a normal projective complex surface with at most rational double points whose minimal resolution is a $\textup{K3}$ surface.

For an even integral lattice $L$, put \[
\operatorname{Roots}(L):=\{v\in L:v^2=-2\}.
\]
An \emph{even overlattice} of an even integral lattice $R$ is an even integral lattice $L$ with $R\subset L\subset R\otimes\mathbb{Q}$, equipped with the extended form, such that the index $[L:R]$ is finite.

If $\rho:\tilde{S}\to S$ is the minimal resolution of a normal projective K3 surface, let $R_S\subset Q_{\tilde{S}}$ be the lattice generated by the exceptional curve classes and let \[
\bar{R}_S:=(R_S\otimes\mathbb Q)\cap H_2(\tilde S;\mathbb Z)
\]
be its primitive closure. 

\begin{theorem}[{\cite[Propositions 2.6 and 2.7]{Shimada-2007}}]
\label{thm:shimada_criterion}
Let $R$ be a negative definite $\operatorname{ADE}$ root lattice of rank $\leq 19$. Then there exists a normal projective $\operatorname{K3}$ surface whose singularities have precisely the Dynkin type of $R$ if and only if there is an even overlattice $L$ of $R$ satisfying \[
\operatorname{Roots}(L)=\operatorname{Roots}(R)
\]
and admitting a primitive embedding into $\Lambda_{\mathrm{K3}}$. Moreover, for every such $L$, the realization may be chosen so that there is an isometry $\bar{R}_S\cong L$ taking $R_S$ onto $R$.
\end{theorem}

We also recall the theorem of Kobayashi--Todorov that does not require projectivity.

\begin{theorem}[{\cite[Corollary 10]{Kobayashi-Todorov-1987}}]
\label{thm:kobayashi_todorov}
Let $S$ be a compact normal complex surface with at most rational double points whose minimal resolution is a $\operatorname{K3}$ surface. Then $e_{\mathrm{orb}}(S)\geq 0$. Equality holds if and only if $S\cong A/G$, where $A$ is a complex $2$-torus and $G$ is a finite group of holomorphic affine transformations acting freely away from finitely many points. In particular, in the equality case $A\to S$ is a finite holomorphic orbifold covering.
\end{theorem}

Here, an affine automorphism of a complex torus $A=\mathbb{C}^2/\Lambda_A$ means an automorphism induced by \[
z\longmapsto Lz+b,\quad
(L\in\mathrm{GL}(2,\mathbb{C}),\, L\Lambda_A=\Lambda_A, \, b\in \mathbb{C}^2).
\]

\section{Proof of Theorem \ref{thm:main}}
\label{sec:vanishing_canonical}

We first treat symplectic $4$-orbifolds with only rational double points and then reduce the general case to this one using Chen's Calabi--Yau cover. We begin by showing that finite covers of the smooth locus extend to orbifold covers with the same properties.

\begin{lemma}\label{lem:rdp_cover_completion}
Let $Y$ be a closed, connected symplectic $4$-orbifold with only rational double points and $K_Y=0\in H^2(Y;\mathbb{Q})$. Then every connected finite symplectic covering of $Y^\circ$ completes to a connected finite symplectic orbifold covering $q:Z\to Y$. Moreover, the orbifold $Z$ is closed, has only rational double points, and satisfies $K_Z=0\in H^2(Z;\mathbb{Q})$.
\end{lemma}
\begin{proof} Near a singular point of $Y$, write the orbifold chart as $B^4/\Gamma$, where $\Gamma\subset\mathrm{SU}(2)$ acts freely away from the origin. A connected component of the restricted covering of $(B^4-\{0\})/\Gamma$ corresponds to a subgroup $N\subset\Gamma$ and is isomorphic, as a covering, to $(B^4-\{0\})/N$. Complete each such component by adding the origin, thereby obtaining the local model $B^4/N$. The resulting map $B^4/N\to B^4/\Gamma$ is an orbifold covering. The pullback of the symplectic form extends over this chart because its uniformizing lift is the original $\Gamma$-invariant symplectic form. These extensions agree with the pullback form on the punctured charts. Since $N\subset\mathrm{SU}(2)$, each newly attached point is either smooth or a rational double point. Finally, the orbifold tangent bundle pulls back under an orbifold covering, so $K_Z=q^*K_Y=0\in H^2(Z;\mathbb{Q})$.
\end{proof}

Next, Bauer's theorem gives the following uniform lower bound for the orbifold Euler characteristic. Although it is weaker than nonnegativity, the bound applies to all finite orbifold covers and it will bound their degrees when the orbifold Euler characteristic is negative.

\begin{lemma}\label{lem:rdp_cy_uniform_lower_bound}
Let $Y$ be a closed, connected symplectic $4$-orbifold with only rational double points and $K_Y=0\in H^2(Y;\mathbb{Q})$. Then \[
e_{\mathrm{orb}}(Y)\geq-\frac{9}{2}.
\]
The same bound holds for every connected finite symplectic orbifold cover of $Y$.
\end{lemma}
\begin{proof}
Let $\rho:M\to Y$ be the minimal symplectic resolution. Then $K_M$ is torsion by Lemma \ref{lem:rdp_canonical_pullback}, so Theorem \ref{thm:Bauer} gives $b_2^+(M)\leq 3$. For each $p\in \operatorname{Sing}(Y)$, let $\ell_p$ be the number of exceptional curves in $\rho^{-1}(p)$ and let $\Gamma_p$ be its isotropy group. Set \[
\delta_p:=\ell_p+1-\frac1{|\Gamma_p|} \ \ \mathrm{and} \ \ L:=\sum_p\ell_p.
\]
Then \[
e_{\mathrm{orb}}(Y)=e(Y)-\sum_{p\in \operatorname{Sing}(Y)} \left(1-\frac{1}{|\Gamma_p|}\right)=e(M)-\sum_{p\in \operatorname{Sing}(Y)}\delta_p.
\]
If $\ell_p=1$, the point is of type $A_1$ and $\delta_p=3/2=3\ell_p/2$. If $\ell_p\geq2$, then $\delta_p<\ell_p+1\leq3\ell_p/2$. Consequently, $\sum_{p\in \operatorname{Sing}(Y)}\delta_p\leq3L/2$. The exceptional lattice is negative definite, so $L\leq b_2^-(M)$. Also, $K_M^2=0$ implies \[
e(M)=\frac32\bigl(b_2^-(M)-b_2^+(M)\bigr).
\]
Combining these observations, we obtain
\[
e_{\mathrm{orb}}(Y) \geq e(M)-\frac{3}{2}L \geq e(M)-\frac{3}{2}b_2^-(M) =-\frac{3}{2}b_2^+(M) \geq-\frac{9}{2},
\]
as desired. The last statement follows from Lemma \ref{lem:rdp_cover_completion}. 
\end{proof}

Next, we record the numerical consequence of Shimada's criterion and the algebraic orbifold BMY inequality that will be used at the end of the covering argument. It will also be used in Section \ref{sec:equality_case}.

\begin{lemma}\label{lem:primitive_ade_complex_realization}
Let $R$ be a negative definite $\operatorname{ADE}$ root lattice of rank $\leq 19$ admitting a primitive embedding into $\Lambda_{\mathrm{K3}}$. Then there exists a normal projective $\operatorname{K3}$ surface $S$ whose singularities have precisely the Dynkin type of $R$, and whose exceptional root lattice is primitive in the intersection lattice of its minimal resolution. In particular, $H_1(S^\circ;\mathbb{Z})=0$ and
\begin{equation}\label{eq:lem:primitive_ade_complex_realization}
24-\sum_p\left(\ell_p+1-\frac1{|\Gamma_p|}\right)\geq0,
\end{equation}
where $p$ ranges over the irreducible $\operatorname{ADE}$ summands of $R$, $\ell_p$ is the rank of the corresponding summand, and $\Gamma_p\subset \mathrm{SU}(2)$ is the isotropy group of that type.
\end{lemma}

\begin{proof}
Apply Theorem \ref{thm:shimada_criterion} with $L=R$. For the resulting minimal resolution $\tilde{S}\to S$, the exceptional lattice $R_S$ is primitive in the unimodular lattice $Q_{\tilde{S}}$. By Lemma \ref{lem:primitive_restriction}, the restriction map $H_2(\tilde{S};\mathbb Z)\to R_S^*$ is surjective. Then the homology long exact sequence
\[
H_2(\tilde{S};\mathbb{Z}) \longrightarrow R_S^* \longrightarrow H_1(S^\circ;\mathbb{Z}) \longrightarrow H_1(\widetilde{S};\mathbb{Z})=0
\]
 gives $H_1(S^\circ;\mathbb{Z})=0$. Finally, the left-hand side of \eqref{eq:lem:primitive_ade_complex_realization} equals $e_{\mathrm{orb}}(S)$ which is nonnegative by Theorem \ref{thm:orbifold_BMY_inequalities}; note that $K_S=0\in H^2(S;\mathbb{Q})$. 
\end{proof}

We can now prove nonnegativity of $e_{\mathrm{orb}}$ in the rational-double-point case. The key idea is to pass to a finite cover of maximal degree, and show that its smooth locus has vanishing first integral homology.

\begin{proposition}\label{prop:rdp_canonical_zero_nonnegativity}
Let $Y$ be a closed symplectic $4$-orbifold with only rational double points and $K_Y=0\in H^2(Y;\mathbb{Q})$. Then $e_{\mathrm{orb}}(Y)\geq0$.
\end{proposition}

\begin{proof}
We may assume that $Y$ is connected. Suppose $e_{\mathrm{orb}}(Y)=-a<0$. For every connected finite symplectic orbifold cover $q:Z\to Y$ of degree $d$, 
Lemma \ref{lem:rdp_cy_uniform_lower_bound} implies \[
-ad=e_{\mathrm{orb}}(Z)\geq-\frac{9}{2},
\]
i.e., $d\leq 9/(2a)$. Choose such a cover $q:Z\to Y$ of maximal degree.
 
We first claim that $H_1(Z^\circ;\mathbb{Z})=0$. Since the first integral homology of the smooth locus $Z^\circ$ is finitely generated, if this group were nonzero, it would have a nontrivial finite cyclic quotient. Then the associated connected finite covering of $Z^\circ$ would complete, by Lemma \ref{lem:rdp_cover_completion}, to a nontrivial connected finite orbifold covering of $Z$. Composing with $q$ would contradict the maximality of its degree. This proves the claim.

Let $\rho: \tilde{Z}\to Z$ be the minimal symplectic resolution and let $R\subset Q_{\tilde{Z}}$ be the exceptional lattice. The homology long exact sequence \[
H_2(\tilde{Z};\mathbb{Z})\longrightarrow R^* \longrightarrow H_1(Z^\circ;\mathbb{Z})=0 \longrightarrow H_1(\tilde{Z};\mathbb{Z})\longrightarrow0
\]
shows that $H_1(\tilde{Z};\mathbb Z)=0$ and that the map $H_2(\tilde{Z};\mathbb Z)\to R^*$ is surjective. Now $H^2(\tilde{Z};\mathbb Z)$ is torsion-free but $K_{\tilde{Z}}\in H^2(\tilde{Z};\mathbb{Z})$ is torsion (Lemma \ref{lem:rdp_canonical_pullback}), so we must have $K_{\tilde{Z}}=0\in H^2(\tilde{Z};\mathbb{Z})$. In particular, $\tilde{Z}$ is spin, which implies that the intersection lattice $Q_{\tilde{Z}}$ is even. 

Writing $(b^+,b^-)=(b_2^+(\tilde{Z}),b_2^-(\tilde{Z}))$, Theorem \ref{thm:Bauer} gives $1\leq b^+\leq3$, while $K_{\tilde{Z}}^2=0$ and $b_1(\tilde{Z})=0$ give \[
0=2e(\tilde{Z})+3\sigma(\tilde{Z})=4+5b^+-b^-, \quad \sigma(\tilde{Z})=b^+-b^-=-4(1+b^+).
\]
Rohlin's theorem forces $\sigma(\tilde{Z})\equiv0\pmod{16}$, so $b^+=3$ and $b^-=19$. Therefore, $e(\tilde{Z})=24$ and $Q_{\tilde{Z}}\cong \Lambda_{\operatorname{K3}}$. Since the restriction $Q_{\tilde{Z}}\to R^*$ is surjective and $Q_{\tilde{Z}}$ is unimodular, Lemma \ref{lem:primitive_restriction} implies that $R$ is primitive in $Q_{\tilde{Z}}$. 

Since $\operatorname{rank}R\leq19$, Lemma \ref{lem:primitive_ade_complex_realization} gives a normal projective K3 surface $S$ with exactly the same singularity type as $Z$. Since the resolutions of both $Z$ and $S$ have Euler characteristic $24$, \[
e_{\mathrm{orb}}(Z) =24-\sum_p\left(\ell_p+1-\frac1{|\Gamma_p|}\right) =e_{\mathrm{orb}}(S)\geq0,
\]
contrary to $e_{\mathrm{orb}}(Z)=-a\deg(q)<0$. This proves the proposition. 
\end{proof}

Finally, it remains to remove the rational-double-point assumption using Chen's Calabi--Yau covering theorem.

\begin{proof}[Proof of Theorem \ref{thm:main}]
We may assume that $X$ is connected. If every singularity of $X$ is a rational double point, we apply Proposition \ref{prop:rdp_canonical_zero_nonnegativity} directly. Otherwise, let $Y\to X$ denote the Calabi--Yau cover of $X$ (Theorem \ref{thm:chen_canonical_cover}) with degree $r=r(X)$. Then Proposition \ref{prop:rdp_canonical_zero_nonnegativity} implies $e_{\mathrm{orb}}(Y)\geq0$, and hence \[
e_{\mathrm{orb}}(X)=\frac1r e_{\mathrm{orb}}(Y)\geq0.
\]
This completes the proof. 
\end{proof}

\section{Proof of Theorem \ref{thm:finiteness}}\label{sec:finiteness}

We begin by recording the Betti-number bounds that will be used throughout. If $M$ is a closed, connected symplectic $4$-manifold with torsion canonical class, then $K_M^2=0$, and hence \[
b_2^-(M)=4-4b_1(M)+5b_2^+(M),\quad
b_2(M)=4-4b_1(M)+6b_2^+(M).
\]
Theorem \ref{thm:Bauer} therefore gives \[
b_2^-(M)\leq19,\quad b_2(M)\leq22.
\]

Next, put \[
B:=\operatorname{lcm}\{m\geq1:\varphi(m)\leq22\},
\]
where $\varphi$ is Euler's totient function. If a prime power $p^a$ divides $m$ with $\varphi(m)\leq22$, then $p^{a-1}(p-1)\leq22$, which leaves only finitely many possibilities for $m$. Therefore, $B$ is a well-defined integer. See also Remark \ref{rmk:choice_of_B} below.

We now prove Theorem \ref{thm:finiteness}. Suppose first that $X$ has only rational double points, including the possibility that $X$ is smooth. Then $r(X)=1$ and its minimal symplectic resolution $\tilde{X}$ has torsion canonical class by Lemma \ref{lem:rdp_canonical_pullback}. Thus $b_2(X)\leq b_2(\tilde{X})\leq22$. Each irreducible summand of the exceptional root lattice has rank at most $b_2^-(\tilde{X})\leq19$, so \eqref{eq:ade_group_orders} gives $|\Gamma_p|\leq120$ at every singular point.

Now suppose $X$ has a singularity that is not a rational double point. Let $q:Y\to X$ be the Calabi--Yau cover of degree $r:=r(X)>1$, and let $M\to Y$ be its equivariant minimal symplectic resolution, as in Theorem~\ref{thm:chen_canonical_cover}. In particular, $K_M=0\in H^2(M;\mathbb{Z})$, and the deck transformations lift to an effective symplectic $\mathbb{Z}_r$-action on $M$. Moreover, either $M$ has the integral homology of a K3 surface or we have $r\in\{3,5\}$. In the latter case, $r$ divides $B$.

In the integral homology K3 case, $e(M)=24$ and $K_M^2=0$ give $\sigma(M)=-16$, and hence $b_2^+(M)=3$. The induced action on $H^2(M;\mathbb{Z})\cong\mathbb{Z}^{22}$ is faithful. Indeed, let $F$ be the kernel of the associated homomorphism $\mathbb{Z}_r\to \operatorname{Aut}(H^2(M;\mathbb{Z}))$. Since $H^1(M;\mathbb{Q})=H^3(M;\mathbb{Q})=0$ and the action preserves orientation, $F$ acts trivially on all rational homology. Theorem \ref{thm:chen_kwasik_rigidity} implies that the action of $F$ on $M$ is trivial, so effectiveness gives $F=\{0\}$.

Choose an integral basis of $H^2(M;\mathbb{Z})$, and let $A\in\mathrm{GL}(22,\mathbb{Z})$ represent a generator of the lifted $\mathbb{Z}_r$-action. Then $A$ has order $r$. If a complex eigenvalue of $A$ has order $m$, then its minimal polynomial over $\mathbb{Q}$, which is the $m$-th cyclotomic polynomial $\Phi_m$, must divide the characteristic polynomial of $A$. Thus, $\varphi(m)=\deg\Phi_m\leq 22$. Since $A$ is diagonalizable over $\mathbb{C}$, its order $r$ is the least common multiple of the orders of its eigenvalues. Hence $r$ divides $B$.

The pullback $q^*:H^2(X;\mathbb{Q})\to H^2(Y;\mathbb{Q})$ is injective by \eqref{eq:cover_pairing}, and resolution adds only negative definite exceptional classes. Consequently, $b_2(X)\leq b_2(Y)\leq b_2(M)\leq22$. For each $p\in\operatorname{Sing}(X)$, the isotropy group at every point over $p$ in $Y$ is $H_p=\Gamma_p\cap\mathrm{SU}(2)$, by Theorem \ref{thm:chen_canonical_cover}\textup{(a)}. If $H_p\neq\{1\}$, the exceptional root lattice over such a point has rank at most $b_2^-(M)\leq19$, so \eqref{eq:ade_group_orders} gives $|H_p|\leq120$. This also holds when $H_p=\{1\}$. Hence
\[
|\Gamma_p|=r_p|H_p|\leq120r\leq120B.
\]

Finally, in both cases, rational Poincar\'e duality applied to the rational homology $4$-manifold $X$ gives \[
e(X)=2-2b_1(X)+b_2(X)\leq24.
\]
Thus, by Theorem \ref{thm:main}, we have \[
0\leq e_{\mathrm{orb}}(X) =e(X)-\sum_{p\in\operatorname{Sing}(X)}\left(1-\frac{1}{|\Gamma_p|}\right) \leq24-\frac{1}{2}|\operatorname{Sing}(X)|.
\]
Hence $|\operatorname{Sing}(X)|\leq48$. Note that there are only finitely many conjugacy classes of subgroups of $\mathrm U(2)$ of order at most $120B$. Since at most $48$ singular points occur, there are only finitely many unordered collections of quotient singularity types, counted with multiplicity. This proves Theorem \ref{thm:finiteness}.

\begin{remark}\label{rmk:choice_of_B}
The prime-power divisibility criterion above gives \[
B=2^5\cdot3^3\cdot5^2\cdot7\cdot11\cdot13\cdot17\cdot19\cdot23
=160{,}626{,}866{,}400.
\]
This common multiple is used only to obtain a uniform bound; no optimality is asserted.
\end{remark}

\begin{proof}[Proof of Corollary \ref{cor:canonical_zero_qhpp}] Put $n:=|\operatorname{Sing}(X)|$. The nonnegativity $e_{\mathrm{orb}}(X)\geq 0$ and the finiteness of the possible singularity types follow from Theorems \ref{thm:main} and \ref{thm:finiteness}. Since $e(X)=3$, the inequality $e_{\mathrm{orb}}(X)\geq 0$ implies that $n\leq6$. If $n=6$, then all six singularities must be of type $A_1$. The minimal symplectic resolution $M$ then has $b_1(M)=0$, $b_2^+(M)=1$, and $b_2^-(M)=6$, giving $K_M^2=2e(M)+3\sigma(M)=3$. This contradicts Lemma \ref{lem:rdp_canonical_pullback}, since $K_X=0$ in rational cohomology. Thus $n\leq5$.
Note that both $e_{\mathrm{orb}}(X)=0$ and $|\operatorname{Sing}(X)|=5$ are attained by contracting a configuration of type $2A_3+3A_1$ on an Enriques surface: the resulting projective $\mathbb{Q}$-homology projective plane has vanishing rational canonical class and $e_{\mathrm{orb}}=0$; see \cite{Jo-Park-2026}.
\end{proof}

\section{Proof of Theorem \ref{thm:zero_euler_smooth_cover}}
\label{sec:equality_case}

We begin with a consequence of Theorem \ref{thm:kobayashi_todorov} for the first homology of the smooth locus.

\begin{lemma}\label{lem:flat_k3_nonzero_h1}
Let $S$ be a normal projective $\operatorname{K3}$ surface. If $e_{\mathrm{orb}}(S)=0$, then $H_1(S^\circ;\mathbb{Z})\neq0$.
\end{lemma}
\begin{proof}
By Theorem \ref{thm:kobayashi_todorov}, we have $S=A/G$, where $A=\mathbb{C}^2/\Lambda_A$ is a complex $2$-torus and $G$ is a finite group of holomorphic affine transformations acting freely away from finitely many points. The group $G$ is nontrivial, since $A$ itself is not a K3 surface. Since the minimal resolution of $S$ is a K3 surface, its nowhere-vanishing holomorphic $2$-form induces a nowhere-vanishing holomorphic $2$-form on $S^\circ$. By Hartogs' theorem \cite{Huybrechts-2005}, the pullback of this $2$-form on $S^\circ$ extends uniquely to a nonzero holomorphic $2$-form on $A$. The extension is $G$-invariant and is a nonzero constant multiple of the translation-invariant form induced by the $2$-form $dz_1\wedge dz_2$ on $\mathbb{C}^2$. Hence every linear part has determinant $1$.

We claim that the abelianization $G^{\mathrm{ab}}$ of $G$ is nontrivial. Let $T$ and $F$ denote the kernel and the image, respectively, of the homomorphism $G\to \mathrm{SL}(2,\mathbb{C})$ sending an affine transformation to its linear part. If $F$ is trivial, then $G=T$, which is abelian since $T$ consists of translations, so we are done. Assume $F$ is nontrivial. It preserves the lattice $\Lambda_A$, and is conjugate to a finite subgroup of $\mathrm{SU}(2)$. In particular, $F$ has no element of order $5$. Indeed, such an element would have eigenvalues $\zeta^a,\zeta^{-a}$, where $\zeta=e^{2\pi i/5}$ and $a\not\equiv0\pmod{5}$, and its real trace would be $4\cos(2\pi a/5)\notin\mathbb{Z}$, which is a contradiction since it must preserve the lattice $\Lambda_A$. Since the binary icosahedral group has an element of order $5$, the classification of finite subgroups of $\mathrm{SU}(2)$ implies that $F$ is isomorphic to a cyclic, binary dihedral, binary tetrahedral, or binary octahedral group. It follows that $F^{\mathrm{ab}}$ is nontrivial, and hence $G^{\mathrm{ab}}$ is also nontrivial, as claimed. Note that the preimage of $S^\circ$ in $A$ is a connected regular covering of $S^\circ$ with deck transformation group $G$. Consequently, there is an epimorphism $\pi_1(S^\circ)\to G$, and this induces an epimorphism $H_1(S^\circ;\mathbb{Z})\to G^{\mathrm{ab}}\neq 0$, completing the proof.
\end{proof}

The realization in Lemma \ref{lem:primitive_ade_complex_realization} has a smooth locus with vanishing first integral homology, so the preceding lemma rules out equality in \eqref{eq:lem:primitive_ade_complex_realization}:
\begin{lemma}\label{lem:primitive_ade_strict}
Let $R$ be a negative definite $\operatorname{ADE}$ root lattice of rank $\leq 19$ admitting a primitive embedding into $\Lambda_{\mathrm{K3}}$. Then \[
24-\sum_p\left(\ell_p+1-\frac{1}{|\Gamma_p|}\right)>0,
\]
with the notation of Lemma \textup{\ref{lem:primitive_ade_complex_realization}}.
\end{lemma}

We now show that, when $e_{\mathrm{orb}}=0$, the degrees of finite orbifold covers are unbounded. We then use the uniform bound on the exceptional rank to prove that a sufficiently large regular cover is smooth.

\begin{proof}[Proof of Theorem \ref{thm:zero_euler_smooth_cover}]
If every singularity of $X$ is a rational double point, set $Y=X$. Otherwise, let $Y\to X$ be the Calabi--Yau cover of Theorem \ref{thm:chen_canonical_cover}. In either case, $Y$ is connected, has only rational double points, and satisfies $K_Y=0\in H^2(Y;\mathbb{Q})$ and $e_{\mathrm{orb}}(Y)=0$. If $Y$ is smooth, we are done, so assume that $\operatorname{Sing}(Y)$ is nonempty.

We claim that the degrees of connected finite orbifold covers of $Y$ are unbounded. Otherwise, choose an orbifold cover $Z\to Y$ of maximal degree. Since $H_1(Z^\circ;\mathbb{Z})$ is finitely generated, any nonzero such group has a nontrivial finite cyclic quotient. The corresponding cover of $Z^\circ$ completes by Lemma \ref{lem:rdp_cover_completion}, contradicting maximality. Hence $H_1(Z^\circ;\mathbb{Z})=0$.

Let $\tilde{Z} \to Z$ be the minimal symplectic resolution and $R\subset Q_{\tilde{Z}}$ its exceptional root lattice. As in the proof of Proposition \ref{prop:rdp_canonical_zero_nonnegativity}, we have $H_1(\tilde{Z};\mathbb{Z})=0$, $K_{\tilde{Z}}=0\in H^2(\tilde{Z};\mathbb{Z})$, $e(\tilde{Z})=24$, $Q_{\tilde{Z}}\cong\Lambda_{\mathrm{K3}}$, and $R$ is primitive in $Q_{\tilde{Z}}$. Since $R$ is negative definite, its rank is $\leq 19$. Lemma \ref{lem:primitive_ade_strict} now implies \[
e_{\mathrm{orb}}(Z)=24-\sum_p\left(\ell_p+1-\frac{1}{|\Gamma_p|}\right)>0,
\]
contrary to $e_{\mathrm{orb}}(Z)=\deg(Z\to Y)e_{\mathrm{orb}}(Y)=0$. This proves the claim.

The degrees of connected finite regular orbifold covers of $Y$ are also unbounded. Indeed, restrict any connected finite orbifold cover to $Y^\circ$, replace its corresponding finite-index subgroup of $\pi_1(Y^\circ)$ by its normal interior, and complete the resulting regular covering using Lemma \ref{lem:rdp_cover_completion}. This produces a regular orbifold cover of at least the original degree.

Put $g:=\max_{p\in\operatorname{Sing}(Y)}|\Gamma_p|$, and choose a connected finite regular orbifold covering $q:W\to Y$ of degree $d>19g$. Let $\tilde{W}\to W$ be its minimal symplectic resolution. Its canonical class is torsion, so $b_2^-(\tilde{W})\leq 19$ as in the beginning of Section \ref{sec:finiteness}. Let $D$ be the deck transformation group of $q$. For $p\in\operatorname{Sing}(Y)$, let $N_p$ denote the kernel of the local monodromy $\Gamma_p\to D$. Every point above $p$ has isotropy group isomorphic to $N_p$, and the number of these points is $d/[\Gamma_p:N_p]$. If $N_p\neq\{1\}$, each of these points contributes at least one exceptional sphere, and the corresponding exceptional configurations span a negative definite
subspace. It follows that \[
b_2^-(\tilde{W}) \geq\frac{d}{[\Gamma_p:N_p]} \geq\frac{d}{|\Gamma_p|}>19,
\]
a contradiction. Thus every $N_p$ is trivial, and $W$ is smooth. Composing $W\to Y$ with $Y\to X$ gives the required covering, with $M=W$. It is clear that $K_M\in H^2(M;\mathbb{Z})$ is torsion and that $e(M)=0$.
\end{proof}

\begin{remark}\label{rmk:positive_euler_no_smooth_cover}
The hypothesis $e_{\mathrm{orb}}(X)=0$ of Theorem \ref{thm:zero_euler_smooth_cover} cannot be removed. To see this, let $S$ be a normal projective K3 surface with exactly one singular point of type $A_1$; such a surface exists by Lemma \ref{lem:primitive_ade_complex_realization}, since $A_1$ admits a primitive embedding into $\Lambda_{\mathrm{K3}}$. Then $K_S=0$ and $e_{\mathrm{orb}}(S)=45/2$. If a connected finite orbifold covering $M\to S$ of degree $d$ had a smooth total space, an orbifold K\"ahler form on $S$ would pull back to a symplectic form on $M$. The canonical class $K_M\in H^2(M;\mathbb{Z})$ would be torsion, and Theorem \ref{thm:Bauer} would give
\[
e(M)=6\bigl(1-b_1(M)+b_2^+(M)\bigr)\leq24.
\]
On the other hand, multiplicativity gives $e(M)=45d/2$, and hence $d\leq16/15<2$. Thus, $d=1$ and $M\cong S$, which is a contradiction since $S$ is not a smooth manifold. Therefore $S$ admits no finite orbifold covering by a smooth manifold.
\end{remark}

% \clearpage
%\bibliography{references}{}

\begin{thebibliography}{BHPV04}

\bibitem[Bau08]{Bauer-2008}
Stefan Bauer.
\newblock Almost complex 4-manifolds with vanishing first {C}hern class.
\newblock {\em J. Differential Geom.}, 79(1):25--32, 2008.

\bibitem[BHPV04]{Barth-Hulek-Peters-VanDeVen-2004}
Wolf~P. Barth, Klaus Hulek, Chris A.~M. Peters, and Antonius {Van de Ven}.
\newblock {\em Compact complex surfaces}, volume~4 of {\em Ergebnisse der Mathematik und ihrer Grenzgebiete. 3. Folge. A Series of Modern Surveys in Mathematics}.
\newblock Springer-Verlag, Berlin, second edition, 2004.

\bibitem[BKR26]{K3-2026}
R.~{\.I}nan{\c{c}} Baykur, Robion~C. Kirby, and Daniel Ruberman, editors.
\newblock {\em {$K3$}: A New Problem List in Low-Dimensional Topology}, volume 295 of {\em Mathematical Surveys and Monographs}.
\newblock American Mathematical Society, Providence, RI, 2026.

\bibitem[Bla95]{Blache-1995}
Raimund Blache.
\newblock The structure of l.c.\ surfaces of {K}odaira dimension zero. {I}.
\newblock {\em J. Algebraic Geom.}, 4(1):137--179, 1995.

\bibitem[CFM08]{Cavalcanti-Fernandez-Munoz}
Gil~R. Cavalcanti, Marisa Fern{\'a}ndez, and Vicente Mu{\~n}oz.
\newblock Symplectic resolutions, {L}efschetz property and formality.
\newblock {\em Adv. Math.}, 218(2):576--599, 2008.

\bibitem[Che20]{Chen-2020}
Weimin Chen.
\newblock On a class of symplectic $4$-orbifolds with vanishing canonical class.
\newblock {\em J. G\"okova Geom. Topol. GGT}, 14:55--90, 2020.

\bibitem[CK07]{Chen-Kwasik-2007}
Weimin Chen and Slawomir Kwasik.
\newblock Symplectic symmetries of $4$-manifolds.
\newblock {\em Topology}, 46(2):103--128, 2007.

\bibitem[ETW26]{Esser-2022}
Louis Esser, Burt Totaro, and Chengxi Wang.
\newblock Calabi-{Y}au varieties of large index.
\newblock {\em Algebr. Geom.}, 13(5):639--666, 2026.

\bibitem[Huy05]{Huybrechts-2005}
Daniel Huybrechts.
\newblock {\em Complex geometry}.
\newblock Universitext. Springer-Verlag, Berlin, 2005.

\bibitem[JP26]{Jo-Park-2026}
Woohyeok Jo and Jongil Park.
\newblock Symplectic $\mathbb{Q}$-homology projective planes with rational double points, 2026.
\newblock arXiv:2611.36352.

\bibitem[JPP26]{Jo-Park-Park-2026}
Woohyeok Jo, Jongil Park, and Kyungbae Park.
\newblock The algebraic {M}ontgomery--{Y}ang problem, 2026.
\newblock arXiv:2607.06686.

\bibitem[KM99]{KeM-1999}
Se\'{a}n Keel and James McKernan.
\newblock Rational curves on quasi-projective surfaces.
\newblock {\em Mem. Amer. Math. Soc.}, 140(669):viii+153, 1999.

\bibitem[KNS89]{KoNS-1989}
Ryoichi Kobayashi, Shu Nakamura, and Fumio Sakai.
\newblock A numerical characterization of ball quotients for normal surfaces with branch loci.
\newblock {\em Proc. Japan Acad. Ser. A Math. Sci.}, 65(7):238--241, 1989.

\bibitem[Kol08]{Kollar-2008}
J\'{a}nos Koll\'{a}r.
\newblock Is there a topological {B}ogomolov-{M}iyaoka-{Y}au inequality?
\newblock {\em Pure Appl. Math. Q.}, 4(2):203--236, 2008.

\bibitem[KT87]{Kobayashi-Todorov-1987}
Ryoichi Kobayashi and Andrey~N. Todorov.
\newblock Polarized period map for generalized {$K3$} surfaces and the moduli of {E}instein metrics.
\newblock {\em Tohoku Math. J. (2)}, 39(3):341--363, 1987.

\bibitem[Kud02]{Kudryavtsev-2002}
S.~A. Kudryavtsev.
\newblock Classification of logarithmic {E}nriques surfaces with {$\delta=2$}.
\newblock {\em Mat. Zametki}, 72(5):715--722, 2002.

\bibitem[Li06]{Li-Quaternionic-2006}
Tian-Jun Li.
\newblock Quaternionic bundles and {B}etti numbers of symplectic 4-manifolds with {K}odaira dimension zero.
\newblock {\em Int. Math. Res. Not.}, pages Art. ID 37385, 28, 2006.

\bibitem[Meg99]{Megyesi-1999}
G.~Megyesi.
\newblock Generalisation of the {B}ogomolov-{M}iyaoka-{Y}au inequality to singular surfaces.
\newblock {\em Proc. London Math. Soc. (3)}, 78(2):241--282, 1999.

\bibitem[Miy77]{Miyaoka-1977}
Yoichi Miyaoka.
\newblock On the {C}hern numbers of surfaces of general type.
\newblock {\em Invent. Math.}, 42:225--237, 1977.

\bibitem[Miy84]{Miyaoka-1984}
Yoichi Miyaoka.
\newblock The maximal number of quotient singularities on surfaces with given numerical invariants.
\newblock {\em Math. Ann.}, 268(2):159--171, 1984.

\bibitem[OZ96]{Oguiso-Zhang-1996}
Keiji Oguiso and De-Qi Zhang.
\newblock On the most algebraic {$K3$} surfaces and the most extremal log {E}nriques surfaces.
\newblock {\em Amer. J. Math.}, 118(6):1277--1297, 1996.

\bibitem[OZ99]{Oguiso-Zhang-1999}
Keiji Oguiso and De-Qi Zhang.
\newblock On the complete classification of extremal log {E}nriques surfaces.
\newblock {\em Math. Z.}, 231(1):23--50, 1999.

\bibitem[Sak80]{Sakai-1980}
Fumio Sakai.
\newblock Semistable curves on algebraic surfaces and logarithmic pluricanonical maps.
\newblock {\em Math. Ann.}, 254(2):89--120, 1980.

\bibitem[Sat57]{Satake-1957}
Ichiro Satake.
\newblock The {G}auss-{B}onnet theorem for {$V$}-manifolds.
\newblock {\em J. Math. Soc. Japan}, 9:464--492, 1957.

\bibitem[Sch23]{Schutt-2023}
Matthias Schütt.
\newblock Moduli of {G}orenstein $\mathbb{Q}$-homology projective planes.
\newblock {\em Journal of the Mathematical Society of Japan}, 75(1):329--366, 2023.

\bibitem[Shi07]{Shimada-2007}
Ichiro Shimada.
\newblock On normal {$K3$} surfaces.
\newblock {\em Michigan Math. J.}, 55(2):395--416, 2007.

\bibitem[Wan11]{Wang-2011}
Fei Wang.
\newblock The classification of log {E}nriques surfaces of rank 18.
\newblock {\em J. Korean Math. Soc.}, 48(4):797--822, 2011.

\bibitem[Yau77]{Yau-1977}
Shing~Tung Yau.
\newblock Calabi's conjecture and some new results in algebraic geometry.
\newblock {\em Proc. Nat. Acad. Sci. U.S.A.}, 74(5):1798--1799, 1977.

\bibitem[Zha91]{Zhang-1991}
De-Qi Zhang.
\newblock Logarithmic {E}nriques surfaces.
\newblock {\em J. Math. Kyoto Univ.}, 31(2):419--466, 1991.

\bibitem[Zha93]{Zhang-1993}
De-Qi Zhang.
\newblock Logarithmic {E}nriques surfaces. {II}.
\newblock {\em J. Math. Kyoto Univ.}, 33(2):357--397, 1993.

\end{thebibliography}
%\bibliographystyle{alpha}

\end{document}